\documentclass[11pt]{article}

\usepackage{amscd,amsmath, amssymb, fancyhdr, epsfig, color, tocloft, mathtools}
\usepackage{dutchcal}

\usepackage[backref=page]{hyperref}
\renewcommand*{\backrefalt}[4]{%
	\ifcase #1 (Not cited.)%
	\or        (Cited on page~#2.)%
	\else      (Cited on pages~#2.)%
	\fi}

\hypersetup{
	colorlinks   = true,
	citecolor    = magenta,
	linkcolor    = blue,
	urlcolor     = magenta	
}

\newcommand{\version}{version 2.0,\ \ Aug. 7, 2024}
\makeatletter
\def\x@arrow{\DOTSB\Relbar}
\def\xlongequalsignfill@{\arrowfill@\x@arrow\Relbar\x@arrow}
\providecommand{\xlongequal}[2][]{%
	\ext@arrow 0099\xlongequalsignfill@{#1}{#2}}
\def\xlongrightarrowfill@{\arrowfill@\relbar\relbar\longrightarrow}
\makeatother

\numberwithin{equation}{section}

\def\eqref#1{(\ref{#1})}

\newcommand{\C}{{\mathbb C}}
\newcommand{\R}{{\mathbb R}}

\newcommand{\6}{\partial}
\def\1{\sqrt{-1}\:}

\newcommand{\cntrct}                
{\hspace{2pt}\raisebox{1pt}{\text{$\lrcorner$}}\hspace{2pt}}
\newcommand{\arrow}{{\:\longrightarrow\:}}

\newcommand{\calo}{{\cal O}}

\renewcommand{\bar}{\overline}
\renewcommand{\phi}{\varphi}
\renewcommand{\epsilon}{\varepsilon}
\renewcommand{\geq}{\geqslant}

\newcommand{\Vol}{\operatorname{Vol}}

\newcommand{\B}{\mathcal{B}}
\newcommand{\SKT}{\mathcal{S\!K\!T}}
\newcommand{\LCK}{\mathcal{L\!C\!K}}

\renewcommand{\Im}{\operatorname{Im}}

\newcounter{Mycounter}[section]
\newcounter{lemma}[section]
\renewcommand{\thelemma}{{Lemma \thesection.\arabic{lemma}}}
\newcommand{\lemma}{%
	\setcounter{lemma}{\value{Mycounter}}
	\refstepcounter{lemma}
	\stepcounter{Mycounter}
	{\noindent \bf \thelemma:\ }}

\newcounter{claim}[section]
\renewcommand{\theclaim}{{Claim \thesection.\arabic{claim}}}
\newcommand{\claim}{%
	\setcounter{claim}{\value{Mycounter}}
	\refstepcounter{claim}
	\stepcounter{Mycounter}
	{\noindent \bf \theclaim:\ }}

\newcounter{sublemma}[section]
\newcounter{corollary}[section]
\renewcommand{\thecorollary}{{Corollary \thesection.\arabic{corollary}}}
\newcommand{\corollary}{%
	\setcounter{corollary}{\value{Mycounter}}
	\refstepcounter{corollary}
	\stepcounter{Mycounter}
	{\noindent \bf \thecorollary:\ }}

\newcounter{theorem}[section]
\renewcommand{\thetheorem}{{Theorem \thesection.\arabic{theorem}}}
\newcommand{\theorem}{%
	\setcounter{theorem}{\value{Mycounter}}
	\refstepcounter{theorem}
	\stepcounter{Mycounter}
	{\noindent \bf \thetheorem:\ }}

\newcounter{conjecture}[section]
\renewcommand{\theconjecture}{{Conjecture \thesection.\arabic{conjecture}}}
\newcommand{\conjecture}{%
	\setcounter{conjecture}{\value{Mycounter}}
	\refstepcounter{conjecture}
	\stepcounter{Mycounter}
	{\noindent \bf \theconjecture:\ }}

\newcounter{proposition}[section]
\newcounter{definition}[section]
\renewcommand{\thedefinition} {{Definition~\thesection.\arabic{definition}}}
\newcommand{\definition}{%
	\setcounter{definition}{\value{Mycounter}}
	\refstepcounter{definition}
	\stepcounter{Mycounter}
	{\noindent \bf \thedefinition:\ }}

\newcounter{example}[section]
\renewcommand{\theexample}{{Example \thesection.\arabic{example}}}
\newcommand{\example}{%
	\setcounter{example}{\value{Mycounter}}
	\refstepcounter{example}
	\stepcounter{Mycounter}
	{\noindent \bf \theexample:\ }}

\newcounter{remark}[section]
\renewcommand{\theremark}{{Remark \thesection.\arabic{remark}}}
\newcommand{\remark}{%
	\setcounter{remark}{\value{Mycounter}}
	\refstepcounter{remark}
	\stepcounter{Mycounter}
	{\noindent \bf \theremark:\ }}

\newcounter{problem}[section]
\newcounter{question}[section]
\renewcommand{\thequestion}{{Question \thesection.\arabic{question}}}
\newcommand{\question}{%
	\setcounter{question}{\value{Mycounter}}
	\refstepcounter{question}
	\stepcounter{Mycounter}
	{\noindent \bf \thequestion:\ }}

\makeatletter

\@addtoreset{equation}{section}
\@addtoreset{footnote}{section}
\makeatother

\def\blacksquare{\hbox{\vrule width 5pt height 5pt depth 0pt}}
\def\endproof{\blacksquare}

\newcommand{\proof}{{\bf Proof: \ }}
\newcommand{\pstep}{{\bf Proof. Step 1: \ }}

\begin{document}
	
	\begin{center}
{\Large\bf  Balanced metrics and Gauduchon cone of locally\\[5mm] conformally K\"ahler manifolds}\\[5mm]
		{\large
			Liviu Ornea\footnote{Liviu Ornea is  partially supported by the PNRR-III-C9-2023-I8 grant CF 149/31.07.2023 Conformal
				Aspects of Geometry and Dynamics.},  
			Misha Verbitsky\footnote{Misha Verbitsky is partially supported by
				FAPERJ 	SEI-260003/000410/2023 and CNPq - Process 310952/2021-2.\\[1mm]
				\noindent{\bf Keywords:} Locally conformally K\"ahler metric, balanced metric, pluriclosed metric, positive current, pseudo-effective cone, Bott--Chern cohomology, Aeppli cohomology, Oeljeklaus--Toma manifold, Kato manifold.\\[1mm]
				\noindent {\bf 2010 Mathematics Subject Classification:} {53C55, 32H04.}
			}\\[4mm]
			
		}
		
	\end{center}

	\hfill
	
	{\small
		\hspace{0.15\linewidth}
		\begin{minipage}[t]{0.7\linewidth}
			{\bf Abstract} \\
A complex Hermitian $n$-manifold $(M,I, \omega)$ is called
locally conformally K\"ahler (LCK) if $d\omega=\theta\wedge\omega$,
where $\theta$ is a closed 1-form, balanced if $\omega^{n-1}$
is closed, and SKT if $dId\omega=0$. We conjecture that any compact 
complex manifold admitting two of these three types
of Hermitian forms (balanced, SKT, LCK) also
admits a K\"ahler metric, and prove partial results
towards this conjecture. We conjecture that
the (1,1)-form $-d(I\theta)$ is Bott--Chern homologous to
a positive (1,1)-current. This conjecture implies that
$(M,I)$ does not admit a balanced Hermitian metric.
We verify this conjecture for all known classes of
LCK manifolds. 
		\end{minipage} 
	}
	\tableofcontents
	

\section{Introduction}


\subsection{Special types of non-K\"ahler Hermitian
  metrics}

Let $(M,I)$ be a connected,  complex manifold of
complex dimension $n\geq 2$. For a Hermitian metric $g$,
we shall denote $\omega(\cdot,\cdot):=g(I\cdot,\cdot)$ its
fundamental two-form. We shall denote by $d^c$ the complex
differential $d^c:=IdI$, where the complex structure acts
on differential forms of degree $k$ by
$I(\eta)(x_1,\ldots,x_k)=(-1)^k\eta(Ix_1,\ldots,Ix_k)$. Sometimes we write $\eta^c$ for $I(\eta)$.

Usually, we work on compact manifolds. 

Further on, we will say ``non-K\"ahler''
referring to a complex manifold not admitting any K\"ahler
structure. A complex manifold admitting a K\"ahler
structure is called ``K\"ahler-type manifold'' if
we do not want to specify the particular K\"ahler metric.

In the realm of non-K\"ahler geometry, the most studied
class of Hermitian metrics is the locally conformally
K\"ahler one. We recall the definition.

\hfill

\definition\label{_LCK_Definition_}
A Hermitian manifold 
$(M,I,g,\omega)$ is called 
{\bf locally conformally K\"ahler} (LCK) if 
there exists a closed 1-form $\theta$ such that $d\omega=\theta\wedge\omega$.
The 1-form $\theta$ is called the {\bf Lee form} 
and its cohomology class {\bf the Lee class}. 

\hfill

\remark 
This definition {is equivalent to
	the existence of a K\"ahler cover $(\tilde
	M,\tilde\omega)\arrow M$ such that the deck group $\Gamma$
	acts on $(\tilde M, \tilde \omega)$ 
	by holomorphic homotheties.} Indeed, suppose that
$\theta$ is exact, $df=\theta$. {Then $e^{-f} \omega$ 
	is a K\"ahler form.} Let $\tilde M$
be a covering such that the pullback $\tilde \theta$ of $\theta$ is exact, 
$df=\tilde \theta$ .
Then the pullback of $\tilde \omega$ is conformal to the K\"ahler form
$e^{-f} \tilde \omega$.

\hfill

An up-to-date account of this kind of geometry appears in  \cite{_OV_book_}.

However, there exist several other interesting  types of non-K\"ahler Hermitian metrics. Below, we recall the definitions of two of these classes, probably the most frequently studied:

\hfill

\definition A Hermitian metric is called {\bf pluriclosed} or {\bf SKT} (shortcut for {\bf strongly K\"ahler with torsion}) if its fundamental two-form satisfies $dd^c\omega=0$ (\cite{_Bismut:connection_}).

\hfill

\remark An Hermitian form $\omega$ is called {\bf
  Gauduchon} if $dd^c(\omega^{n-1})=0$, where $n=\dim_\C M$.
In complex dimension 2, the SKT condition is
equivalent to the Gauduchon condition, and
hence SKT metrics  exist in any conformal class of metrics
on a compact complex surface. In higher dimension, a
general existence result is not known, but many examples
are found on certain nilmanifolds (see, for example, the
papers of Anna Fino and collaborators). Moreover SKT
metrics are connected with certain theories in theoretical
phyisics and are essential in the literature about
generalized complex geometry and generalized K\"ahler
structures (see the works of Cavalcanti, Gualtieri,
Hitchin etc.).

\hfill

\remark Recall that a two-form $\omega$ is called {\bf taming} or 
{\bf symplectic-Hermitian} if it is 
the (1,1)-part of a symplectic form. 
Clearly, a symplectic-Hermitian form is pluriclosed.
The converse is false. Indeed, {there are no examples of 
	symplectic-Hermitian form on non-K\"ahler compact
	complex manifolds}; it is conjectured they don't exist (\cite{_Streets_Tian_}).

\hfill

\definition A Hermitian metric is called {\bf balanced} or {\bf semi-K\"ahler} if its fundamental two-form satisfies $d\omega^{n-1}=0$, $n=\dim_\C M$. Equivalently, the metric is balanced if its two-form is coclosed (\cite{_Michelson_}).

\hfill

\remark In complex dimension 2, balanced and K\"ahler
metrics are equivalent. However, in any complex dimension
greater than 3 there exist compact complex manifolds
carrying balanced metrics, but no K\"ahler
metric. Balanced metrics exist on twistor spaces of
4-dimensional ASD Riemannian manifolds,
\cite{_Michelson_}, while only the twistor spaces of $\C
P^2$ and $S^4$ are K\"ahler (\cite{_Hitchin_}). Moreover,
all Moishezon (and, more generally, Fujiki class C)
manifolds are balanced, because the existence of
balanced metrics is preserved by bimeromorphisms
(\cite{_Alessandrini_Bassanelli:bimero_}).

\hfill

Denote by  $\B$, resp. $\SKT$, resp. $\LCK$ the class of balanced, resp. SKT, resp. LCK manifolds. A natural question arises:

\hfill

\question\label{_different_existence_} Is it possible that
a compact complex non-K\"ahler manifold $(M,I)$ admit   metrics in
two of the above classes? We recall that ``non-K\"ahler manifold''
means ``manifold not admitting a K\"ahler metric''.

\hfill 

In general, the answer is negative, at least for a smaller
class of LCK metrics, namely the Vaisman ones (see
below), for example \cite{_Angella_Otiman_}. More precisely, one proves that the simultaneous
existence of metrics in two different classes forces the
manifold to be of K\"ahler type. Most computations are carried out
on nilmanifolds, see
e.g. \cite{_Fino_Vezzoni:SKT_,_Fino_Vezzoni:correction_}
(for metrics in $\B$ and $\SKT$),
\cite{_Ornea_Otiman_Stanciu_} (for metrics in $\B$ and
$\LCK$), or \cite{_Arroyo_Nicolini_} for the general existence of invariant SKT metrics on $k$-step nilmanifolds.  We dare to conjecture:

\hfill

\conjecture
Let $M$ be a compact complex manifold which admits
Hermitian forms $\omega_1$ and $\omega_2$ which belong to 
two classes in the set $\{\B,\LCK,\SKT\}$. 
Then $M$ admits a K\"ahler metric.

%

\hfill

The aim of this paper is to give a partial negative answer to \ref{_different_existence_}. Namely, we shall prove that balanced metrics cannot exist on any of the known examples of compact LCK manifolds (\ref{_no_balanced_on_lck_Corollary}).

\subsection{Positive currents in LCK geometry}

For an introduction to positive currents
we refer to \cite{_Demailly:CADG_}.

The notion of ``K\"ahler rank'' was introduced first
in \cite{_HL:intrinsic_}, and further developed
in \cite{_Brunella:Inoue_,_Chiose_Toma:Kahler_rank_}.
This notion is a powerful tool used to classify
non-K\"ahler complex surfaces. A non-K\"ahler complex 
surface $M$ is said to be of K\"ahler rank 1 
if it admits a smooth, positive form (which will
a posteriori be semi-positive of rank 1).
Otherwise, it is said that $M$ has K\"ahler rank 0.

The Bott--Chern cohomology group $H^{1,1}_{BC}(M)$
of a complex manifold is the space of closed
(1,1)-forms (or currents) up to $dd^c(C^\infty M)$.
It is well known (see e.g. \cite{_ovv:surf_})
that the kernel of the tautological map $\tau:\; H^{1,1}_{BC}(M)\arrow H^2(M,\R)$
is 1-dimensional for any non-K\"ahler surface.
If $M$ is a non-K\"ahler complex surface,
the kernel of $\tau$ can be represented
by a positive closed (1,1)-current 
(\cite[Th\'eor\`eme 6.1]{_Lamari_}).
This current can be chosen as a smooth positive (1,1)-form
for surfaces of K\"ahler rank 1, and it 
is always singular for surfaces of K\"ahler rank 0
(\cite{_Brunella:Inoue_,_Chiose_Toma:Kahler_rank_}).

In \cite{_ovv:surf_}, it was shown that
all non-K\"ahler surfaces admit an LCK structure,
with the exception of one of the three Inoue surfaces and 
possible non-Kato class VII surfaces
(which conjecturally don't exist).
The proof is based on Lamari's Theorem 6.1;
the positive current constructed by Lamari
is Bott--Chern homologous to $d^c\theta$, where $\theta$
is the Lee form of the LCK structure on
the surface.

On Vaisman and OT-manifolds, equipped with the
most natural LCK metric, the form $d^c\theta$
is actually semi-positive (\cite{_Va_torino_,_OV:OT_};
see also equation \eqref{_dc_theta_}).

This observation leads us to 
\ref{_pseudo_ef_conjecture_} 
which claims essentially that $d^c\theta$
is Bott--Chern homologous to a positive current
on any LCK manifold.

We verify this conjecture for all known
examples of LCK manifolds
(\ref{_pseff_conj_true_for_known_classes_Theorem_}).

This conjecture easily implies that a strict LCK manifold
does not admit a balanced metric 
(\ref{_ps_eff_conj_then_not_balanced_Theorem_}).

\section{A metric which belongs to two of
the classes $\B$, $\LCK$, $\SKT$ is K\"ahler} 

By using only the very definitions of balanced, LCK and SKT metrics, one can prove that on a compact complex manifold, a Hermitian metric cannot belong to the intersections $\B\cap\LCK$, $\B\cap\SKT$, or $\LCK\cap\SKT$ unless it is K\"ahler.

The following \ref{_pairs_lck_skt_b_imply_K:Theorem_} is
partly known, but we provide a unified proof. For other
proofs showing that a Hermitian metric  which is both
balanced and  SKT on a compact complex manifold is
K\"ahler, see \cite[Theorem 1.3]{_Ivanov_Papadopoulos_},
\cite[Proposition 2.6]{_Dinew_Popovici_}. See also
\cite{_Chiose_Rasdeaconu_Suvaina_} for the relation
between SKT and balanced metrics. 
Also, in \cite[Remark 1]{_Alexandrov_Ivanov_} it
was shown that a non-K\"ahler metric on a compact complex
manifold of dimension $>2$ cannot
be both SKT and LCK.

\hfill

\theorem\label{_pairs_lck_skt_b_imply_K:Theorem_}
Let  $(M,I,\omega)$ be a compact Hermitian $n$-manifold, $n\geq 3$.
Assume that $\omega$ belongs to two of the classes $\{\B,\LCK,\SKT\}$. 
 Then $\omega$ is K\"ahler.
 
 \hfill
 
\pstep {\bf Assume  $(M,I,\omega)$ is SKT and LCK.} 
Let $\theta\in \Lambda^1(M)$ be the Lee form of $\omega$ as an LCK metric:  $d\omega=\theta\wedge\omega$.
Let $\theta^c:= I(\theta)$.
Then $d^c \omega = I^{-1} d I(\omega)=I^{-1} (\theta \wedge \omega)=
-\theta^c\wedge \omega$. This gives
\begin{equation}\label{_dc_d_omega_}
	 0=d^cd \omega = d^c (\theta \wedge \omega)=
d^c\theta \wedge \omega - \theta \wedge d^c \omega
= d^c\theta \wedge \omega + \theta \wedge \theta^c\wedge  \omega 
\end{equation}
Since $\dim_\C M > 2$, the multiplication map
$\eta \mapsto \eta \wedge \omega$ is injective, hence
\eqref{_dc_d_omega_} implies that $d^c\theta=-\theta \wedge \theta^c$.
Then 
\begin{equation*}
\begin{split}		
dd^c \omega^{n-1} &= (n-1) d^c\theta\wedge \omega^{n-1}- 
(n-1)^2\theta \wedge \theta^c\wedge \omega^{n-1}\\ 
&=-(n-1) (n-2) \theta \wedge \theta^c\wedge \omega^{n-1}.
\end{split}
\end{equation*}
However, $\theta \wedge \theta^c\wedge \omega^{n-1}= 2n |\theta|^2 \omega^n$.
This brings a contradiction:
\begin{equation*}
	\begin{split}
0 &= \int_M dd^c \omega^{n-1} =
\int_M -(n-1) (n-2) \theta \wedge \theta^c\wedge \omega^{n-1}\\
&= - \frac{(n-1) (n-2)}2n \int_M |\theta|^2\wedge \omega^{n}
	\end{split}
\end{equation*}
The last integral vanishes if and only if $\theta=0$, hence
$\omega$ is closed.
 
\bigskip

{\bf Step 2: 
		Assume  $(M,I,\omega)$ is balanced and LCK.}
	Then $0=d\omega^{n-1} = (n-1)\theta \wedge \omega^{n-1}$.
	However, the multiplication map
	$\eta \mapsto \eta \wedge \omega^{n-1}$ is an isomorphism for all $\eta$
	and any Hermitian $\omega$, hence again
	$\theta=0$ (see also \cite[Exercise 4.16]{_OV_book_}).

\bigskip

{\bf Step 3: 
	Assume  $(M,I,\omega)$ is balanced and SKT.} Then
$d(\omega^{n-1})=(n-1) d\omega \wedge \omega^{n-2}=0$. Moreover,
	$dd^c \omega=0$, hence $d\omega$ and $d^c\omega$
	are $d$ and $d^c$-closed. The equation 
	$d\omega \wedge \omega^{n-2}=0$ implies that $d\omega$ is 
	{\bf primitive}, that is, it satisfies $\Lambda_\omega(d\omega)=0$,
	where $\Lambda_\omega= L_\omega^*$, and $L_\omega(\eta) := \omega\wedge \eta$. 
	This form is of Hodge type $(1,2)+(2,1)$  because $\omega$ is of type (1,1),
	and the de Rham differential shifts the Hodge grading at most by 1.
	
	Recall that by the Hodge-Riemann relations,
	any primitive $(1,2)+(2,1)$ real form $\alpha$
	satisfies $\alpha \wedge I(\alpha)\wedge \omega^{n-3}= -C|\alpha|^2\omega^n$,
	where $C$ is a positive rational constant.
	
	Let $\alpha:= d\omega$.  
	Since $\omega$ is SKT, we have
\begin{equation}\label{_ddc_omega_n-1_}
	\begin{split}
		0&=dd^c(\omega^{n-1})=  (n-1)(n-2) d\omega\wedge d^c\omega\wedge \omega^{n-3}\\
	&=-(n-1)(n-2) C|\alpha|^2\omega^n. 
	\end{split}
\end{equation}
This follows by the previous remark, since $\alpha$ is a
primitive $(1,2)+(2,1)$-form (see also \cite[Exercise 4.24]{_OV_book_}).  It follows again that $d\omega=0$. \endproof

\section{LCK manifolds}

Clearly, the LCK condition (\ref{_LCK_Definition_})
is conformally invariant, and hence conformally equivalent
LCK metrics give rise to cohomologous Lee forms. We are
interested in compact LCK manifolds with non-exact Lee
form, called {\bf strict LCK}. Usually, this condition
is tacitly assumed.
 The following result states the dichotomy between these manifolds and the K\"ahler ones:

\hfill

\theorem {\rm (Vaisman, \cite{_Vaisman_trans_})} 
A compact LCK manifold $(M,I, \theta)$ with non-exact Lee form does not
		admit any K\"ahler structure.

\hfill

\proof
Suppose that $(M,I)$ admits a K\"ahler structure.
By $dd^c$-lemma, $d^c\theta=dd^cf$ for some function $f$.
This is impossible, as the following lemma implies.

\hfill

\lemma\label{_dd^c_Lee_Lemma_}
Let $(M,\omega, \theta)$ be a compact LCK manifold.
Suppose that $d^c \theta$ is $dd^c$-exact. Then $\theta$ is exact.

\hfill

\proof
Suppose that $d^c\theta= d^c df$.
Replacing $\omega$ by $\omega_1:=e^{-f}\omega$ gives an LCK structure
with LCK form $\theta_1=\theta - df$. Then $\theta_1$ is closed
and $d^c$-closed, giving 
$dd^c (\omega_1^{n-1})= (n-1)^2 \theta_1 \wedge \theta_1^c \wedge \omega_1^{n-1}$
The latter form is a volume form, strictly positive at all points
where $\theta_1\neq 0$, hence it cannot be exact unless
$\theta_1=0$.
\endproof

\hfill

Next, we present the main constructions used
to produce LCK manifolds.

\subsection{Blow-up of an LCK manifold}

We owe to Tricerri and Vuletescu the following fundamental result:

\hfill

\theorem{\rm (\cite{_Tricerri_}, \cite{_Vuletescu_})} The blow-up at points preserves the LCK class.

\hfill

\remark However, the blow-up along subvarieties does not preserve the LCK class unless the Lee form is exact along the submanifold, \cite{_OVV_blow_up_}.

\subsection{LCK manifolds with potential}

\definition
An LCK manifold $M$ is called {\bf  LCK manifold with LCK potential}
if the K\"ahler form on the universal cover $\tilde \omega$ of 
$\tilde M$ has a positive K\"ahler potential 
$\phi:\tilde M\to\R^{>0}$ such that the action of
$\pi_1(M)$ multiplies this function by a positive constant.

\hfill 

\remark \label{_defo:remark_} A small deformation of an LCK manifold
might be non-LCK, \cite{_Belgun_}. However, a small
	deformation of LCK with potential is LCK with potential, \cite{_OV_lckpot_}.

\hfill 

\example 
All Hopf manifolds, linear and non-linear, admit an LCK structure
	with LCK potential, \cite{_OV_lckpot_,_OV_non_linear_}.
	
\hfill

\theorem {\rm (\cite[Theorem 13.22]{_OV_book_})\label{_pot_iff_subma_in_Hopf_}}
A compact complex manifold $M$, $\dim_\C M> 2$
admits an LCK potential if and only if
	$M$ admits a holomorphic embedding
	to a Hopf manifold.
	
\hfill

\remark This property
{can be used instead of the definition.}

\hfill

\remark In complex dimension 2, \ref{_pot_iff_subma_in_Hopf_} is still true if we assume the GSS conjecture, \cite[Theorem 25.43]{_OV_book_}.

\subsection{Vaisman manifolds}

In this subsection we recite a few classical results about
Vaisman manifolds. For the reference see
\cite{_OV_book_}.

The best understood (and known) subclass of LCK manifolds with potential is formed by those having the Lee form parallel with respect to the Levi-Civita connection of the LCK metric. Indeed, if $\pi:\ (\tilde M,\tilde\omega)\to (M,\theta,\omega)$ is the universal cover of such an LCK manifold, then one easily verifies that the $\tilde\omega$-squared norm of $\pi^*\theta$ is an LCK potential for $\tilde\omega$. 

\hfill

\example Almost all non-K\"ahler compact complex surfaces are LCK,
see e.g. \cite{_ovv:surf_}. Diagonal Hopf surfaces and
Kodaira surfaces are Vaisman, but Kato surfaces and Inoue
surfaces, as well as OT manifolds (see below) are not Vaisman. In any dimension, all diagonal
Hopf manifolds are Vaisman, see \cite{_OV_pams_}, but non-diagonal Hopf manifolds are not Vaisman. Also, all compact complex submanifolds of a Vaisman manifold are Vaisman (\cite[Corollary 7.35]{_OV_book_}).

\hfill

On a Vaisman manifold, the Lee form can be assumed of norm
1. Then one can prove the following fundamental identity 
(\cite[Formula (2.8)]{_Va_torino_}, \cite[Proposition 8.3]{_OV_book_}):
\begin{equation}\label{_dc_theta_}
	\omega=d^c\theta+\theta\wedge\theta^c.
\end{equation}

\remark Since a non-diagonal Hopf manifold is a small deformation of a diagonal Hopf manifold, which is Vaisman, it follows that small deformations do not preserve the Vaisman property. However, due to \ref{_defo:remark_}, a small deformation of a Vaisman manifold is still LCK (with potential). 

\subsection{Oeljeklaus--Toma manifolds}

The construction of Oeljeklaus--Toma manifolds,
sometimes abbreviated as OT manifolds, is related to number
theory. Their construction generalizes the
construction of the Inoue surfaces $S^0$ in class VII. We
refer to the monograph \cite{_Neukirch_} for the necessary
background in number theory.

Let $K$ be a number field which has $2t$
complex embeddings denoted $\tau_i, \bar \tau_i$ and 
$s$ real ones denoted $\sigma_i$,  $s>0$, $t>0$.

Let $\calo_K^{*,+}:= \calo_K^*\cap \bigcap_i \sigma^{-1}_i(\R^{>0})$.
Choose in $\calo_K^{*,+}$
a free abelian subgroup $\calo_K^{*,U}$ of rank
$s$ such that the quotient
$\R^s/\calo_K^{*,U}$ is compact, where
$\calo_K^{*,U}$ is mapped to $\R^s$ as
$\xi \arrow \big(\log(\sigma_1(\xi)), ..., \log(\sigma_s(\xi))\big).$ The existence of such a subgroup is proven in \cite{_Oeljeklaus_Toma_}.
Let $\Gamma:= \calo^+_K\rtimes \calo_K^{*,U}$, and $H^s$ be the upper half-plane $\{(z_1,\ldots,z_s)\ ;\ \Im(z_1)>0,\ldots,\Im(z_s)>0\}$.

\hfill

\definition {\rm (\cite{_Oeljeklaus_Toma_})} 
An {\bf  Oeljeklaus--Toma (OT) manifold} is a quotient
$\C^t \times H^s/\Gamma$, where
$\calo^+_K$ acts on $\C^t \times H^s$
as 
\[ \zeta(x_1,..., x_t, y_1, ..., y_s) = \bigg(x_1+
\tau_1(\zeta), ..., x_t + \tau_t(\zeta), 
y_1+\sigma_1(\zeta), ..., y_s+\sigma_s(\zeta)\bigg),
\]
and $\calo_K^{*,U}$ as 
\[
\xi(x_1,..., x_t, y_1, ..., y_s)=
\bigg (x_1,..., x_t,\sigma_1(\xi) y_1,  ..., \sigma_s(\xi) y_s\bigg)
\]

\hfill

\theorem {\rm (\cite{_Oeljeklaus_Toma_})}  
The OT-manifold  $M:=\C^t \times H^s/\Gamma$
{is a non-K\"ahler compact complex manifold,} without any non-constant
meromorphic functions. 

\hfill

\theorem An Oeljeklaus--Toma manifold admits an LCK metric if and only if $t=1$ (\cite{_Oeljeklaus_Toma_} for the sufficiency of the condition, \cite{_Deaconu_Vuletescu_} for the necessity).

\subsection{Kato manifolds}

\definition
Let $B$ be an open ball in $\C^n$, $n>1$, and
$\tilde B \stackrel \pi\arrow B$ be a bimeromorphic, 
holomorphic map, which is an isomorphism outside
of a compact subset. Remove a small ball in $\tilde B$
and glue it to the boundary of $\tilde B$, extending
the complex structure smoothly (and holomorphically)
on the resulting manifold, denoted by $M$. 
Then $M$ is called {\bf a Kato manifold.}

\hfill

\theorem {\rm (\cite{_Brunella:Kato_})}
Suppose that $M$ is a Kato manifold obtained from 
$\tilde B \stackrel \pi\arrow B$ with $\tilde B$ K\"ahler.
{Then $M$ is LCK.}

\hfill

Also, Kato manifolds can be obtained as 
limits of modifications of Hopf manifolds:

\hfill

\theorem {\rm (\cite{_Kato:announce_}, \cite{_Dloussky:Kato_})}
Let $M$ be a Kato manifold. Then there exists
a family $M_t$ of complex manifolds over a punctured disk
{such that $M= M_0$ and all other $M_t$ are bimeromorphic to
	a Hopf manifold.}

\hfill

\definition
Let $M$ be a complex manifold, and $\Gamma\subset M$ be an open
subset which is isomorphic as a complex manifold to a small
neighbourhood of a sphere $S^{2n-1} \subset \C^n$. The set $\Gamma$
is called {\bf a global spherical shell} if the complement
$M \backslash \Gamma$ is connected.

\hfill

\theorem {\rm (\cite{_Kato:announce_})}
Let $M$ be a compact complex manifold. Then
{$M$ is a Kato manifold if and only if it contains a 
	global spherical  shell.}

\hfill

More recently, the LCK structure of Kato manifolds was
investigated in \cite{_IOP:Kato_,_IOPR:toric_Kato_}.

\section{The Lee cone versus the Gauduchon cone}

\subsection{The Gauduchon cone}

\definition
A Hermitian metric $\omega$ on a complex $n$-manifold is 
called {a \bf Gauduchon metric} if its fundamental form satisfies  $dd^c(\omega^{n-1})=0$.

\hfill

One of { the very few statements}
which is valid (and useful) to all compact
complex manifolds is the following theorem of Gauduchon:

\hfill 

\theorem {(\cite{_Gauduchon_})} In any conformal class of Hermitian metrics on a compact complex manifold there exists a Gauduchon metric, unique up to a constant multiplier.

\hfill

\definition
Let $M$ be a complex manifold, and $\omega$ a Hermitian form of a Gauduchon metric.
{\bf The Gauduchon form} of $M$ is $\omega^{n-1}$.

\hfill

\claim\label{_Gauduchon_form_is_power_Claim_}
Fix a positive volume form $\Vol$ on $M$.
A form $\eta \in \Lambda^{n-1, n-1}(M,\R)$
defines a Hermitian form on $\Lambda^1(M)$
taking $x, y\in T^*M$ to $\frac{\eta\wedge x \wedge y}{\Vol}$.
{Then this Hermitian form on $\Lambda^1 M$ is positive definite if
	and only if $\eta = \alpha^{n-1}$,} where $\alpha$ is a
Hermitian form on $TM$.

\hfill

\remark From this claim it follows that the set of all Gauduchon forms
	is a convex cone in $\Lambda^{n-1, n-1}(M,\R)$.
Indeed, a sum of closed, strictly positive $(n-1,n-1)$-forms is the Gauduchon
form of a metric, as \ref{_Gauduchon_form_is_power_Claim_} implies.

\hfill

Recall that  the Aeppli cohomology group $H^{p,q}_{AE}(M)$ 
is the space of $dd^c$-closed $(p,q)$-forms modulo
$\6(\Lambda^{p-1,q} M)+ \bar\6(\Lambda^{p,q-1} M)$. 
The Aeppli cohomology is finite-dimensional for $M$ compact. Moreover, 
the natural pairing 
$H^{p,q}_{BC}(M) \times H^{n-p,n-q}_{AE}(M)\arrow
H^{2n}(M)=\C$, taking $x, y$ to $\int_M x\wedge y$
is non-degenerate and identifies $H^{p,q}_{BC}(M)$
with the dual $H^{n-p,n-q}_{AE}(M)^*$ (\cite{_Angella_}). Now we can give:

\hfill

\definition
{\bf The Gauduchon cone} of a compact complex
$n$-manifold is the set of all classes
$\omega^{n-1}\in H^{n-1, n-1}_{AE}(M)$ of all Gauduchon
forms.

\hfill

\definition
The {\bf  pseudo-effective cone}
$P\subset  H^{1,1}_{BC}(M)$
is the set of Bott--Chern classes of all
positive, closed (1,1)-currents.

\hfill

\theorem  \label{_Lamari_pseff_Theorem_}
On a compact complex manifold, the Gauduchon cone is dual to the pseudo-effective cone.\\
\proof \cite[Lemme 3.3]{_Lamari_}; see also \cite{_Popovici_Ugarte_}. \endproof

\subsection{The Lee cone of an LCK manifold}

\definition
Let $M$ be a compact complex manifold 
admitting an LCK structure. {\bf The Lee cone} of $M$ 
is the set of all classes $[\theta] \in H^1(M, \R)$
such that $\theta$ is the Lee form for an LCK structure on $M$.

\hfill

\remark The ``Lee cone'' might be a bad term, because generally speaking the Lee cone is not a cone.
This set can be even discrete, for example on Inoue surfaces $S^0$  (\cite{_Otiman:Inoue_}). 
However, for Vaisman manifolds and, more generally, 
LCK manifolds with potential, the Lee cone is half-space, and {\em ipso facto} a cone
(\cite{_OV_tohoku_}). 

\hfill

On compact LCK manifolds with potential, one can characterize the Lee cone using the {\em degree} of a cohomology class in $H^1(M,\R)$ with respect to a fixed Gauduchon form as follows.

Let $M$ be a compact complex manifold, and $\omega$ a 
Gauduchon form. Consider the natural map 
$$H^1(M, \R) \ni [\alpha] \arrow H^{1,1}_{BC}(M)$$
which takes the class of a closed real 1-form $\alpha$ to the Bott--Chern
class of $d^c\alpha$. This makes sense because locally, $\alpha=df$, hence
$d^c\alpha$ is a (1,1)-form. 

\hfill

\definition The degree of $[\alpha]$ is $\deg_\omega \alpha:=
	\int_M\omega^{n-1} \wedge d^c\alpha$.
	
\hfill

\remark Recall that $dd^c=-d^cd$. Then $\deg\alpha=\deg(\alpha+dh)$ and hence $\deg_\omega$ is a well defined functional on the first cohomology group $H^1(M, \R)$. Moreover, observe that $d^c\alpha\in \Lambda^2(M)$ is always exact.

\hfill

\theorem\cite[Theorem 8.4]{_OV_tohoku_}
Let $M$ be a compact LCK manifold with potential,
and $u$ a Gauduchon form. Then its Lee cone
	is the set of all $\alpha \in H^1(M,\R)$ such that
	$\deg_u \alpha>0$. 
	
\hfill

We conjecture that half of this characterization is true on all LCK manifolds:

\hfill

\conjecture\label{_pseudo_ef_conjecture_} 
Let $\theta$ be a Lee class on a compact LCK manifold, 
and $u$ a Gauduchon metric. Then $\deg_u \theta >0$.

\hfil

\remark
This conjecture is non-trivial, because the choice of the Gauduchon
metric $u$ is arbitrary. When $\omega$ is an LCK form such that $d\omega =\omega\wedge\theta$,
it is easy to check that $\deg_\omega \theta >0$. Indeed, 
\[ 
0 =dd^c(\omega^{n-1})= (n-1)d(\omega^{n-1}\wedge \theta^c)=
(n-1)^2 \omega^{n-1}\wedge \theta\wedge \theta^c + (n-1) \omega^{n-1} d\theta^c.
\]
On the other hand, $d\theta ^c =- d^c\theta$ because $\theta=df$ locally,
and $dd^c=-d^c d$. Then $(n-1)\deg_\omega \theta= (n-1)^2 \int \omega^{n-1}\wedge \theta\wedge \theta^c$,
and $\int \omega^{n-1}\wedge \theta\wedge \theta^c$ is the $L^2$-norm of $\theta$.

\hfill

\remark Recall that an element of $H^{1,1}_{BC}(M)$ is {\bf pseudo-effective} if it is the class of a positive closed (1,1)-current. By \ref{_Lamari_pseff_Theorem_},
the pseudo-effective cone is dual to the Gauduchon cone, that is, a Bott--Chern class
is pseudo-effective if and only if it evaluates positively on all Gauduchon forms. 
 Then \ref{_pseudo_ef_conjecture_} is equivalent to the pseudo-effectivity of the
Bott--Chern class of $d^c\theta$.

\hfill

We still don't know how to prove \ref{_pseudo_ef_conjecture_}. However, we can prove a weaker version:

\hfill

\theorem\label{_pseff_conj_true_for_known_classes_Theorem_}
Let $(M,\omega, \theta)$ be a compact, non-K\"ahler LCK-manifold of complex dimension greater than 3, which 
is bimeromorphic to either  LCK with potential,  Oeljeklaus--Toma or Kato
manifold, and  $u$ a Gauduchon metric. Then 
	$\deg_u \theta > 0$.

\hfill

\pstep
The statement is true for free for OT and Vaisman manifolds,
because $d^c\theta$ is a positive, exact (1,1)-form.
For Vaisman manifolds $d^c\theta$ is positive as shown
in \cite{_Verbitsky:Vanishing_LCHK_}. For OT manifolds
this form is positive as shown in \cite{_OV:OT_}.
It is also true for all manifolds which are bimeromorphic
to Vaisman and OT manifolds, because these manifolds are always 
obtained as blow-ups of Vaisman and OT manifolds, \cite{_OV_bimero_}, and
the pullback of a positive form is positive.
This implies that $\deg_u \theta >0$
in all these situations.

\hfill

{\bf Step 2:}
Suppose we have a smooth family $(M_t,\omega_t, \theta_t)$ of LCK manifolds
such that the statement of the theorem is true for all
$t\neq 0$. Fix a Gauduchon metric $u_0$ on the central
fiber $M_0$. We can extend it $u_0$ to a smooth family
of Gauduchon metrics using the Gauduchon theorem. Indeed,
an Hermitian form $u_0$ in the central fiber $M_0$ can be extended to
a family $u_t'$ of Hermitian forms on the fibers in the 
neighbourhood of $M_0$. Applying Gauduchon theorem, we find a
Gauduchon form $u_t$ conformal to $u'_t$.
Since the solution of the Gauduchon equation
depends smoothly on the data, the form $u_t$
smoothly depends on $t$.

On all fibers except the central, we have
	$\deg_{u_t} \theta_t \geq 0$, hence the same is true on
	the central fiber.
	
\hfill

{\bf Step 3:} Let $M_0$ be 
an LCK manifold with potential. By \cite[Theorem 16.21]{_OV_book_},
$M_t$ has a deformation $(M_t,\omega_t, \theta_t)$ 
such that $(M_t,\omega_t, \theta_t)$ is Vaisman
for all $t\neq 0$. The Kato manifold has a deformation
$(M_t,\omega_t, \theta_t)$ 
with all the fibers except the central one biholomorphic to
a blown-up Hopf (\cite{_Kato:announce_}). However, a blown-up Hopf
manifold is bimeromorphic to an LCK manifold with potential,
hence it also satisfies the statement of the theorem.
Applying Step 2, we obtain that $\deg_u \theta\geq 0$
for all Gauduchon metrics $u$.

\hfill

{\bf Step 4:} It remains to prove that the inequality
$\deg_u \theta\geq 0$ is strict. \ref{_Lamari_pseff_Theorem_},
together with Step 3, implies that $d^c \theta$ is pseudo-effective,
that is, can be represented by a positive current $\Theta$.
If the integral $\int_M \Theta\wedge u^{n-1}$ vanishes
for some Hermitian form $u$, this implies that $\Theta=0$,
and hence the Bott--Chern class of $d^c\theta$ vanishes,
Then $d^c\theta=dd^c f$, for some $f\in C^\infty M$.
This is impossible by \ref{_dd^c_Lee_Lemma_}.
\endproof

\section{Non-existence of balanced metrics on LCK manifolds}

\theorem\label{_ps_eff_conj_then_not_balanced_Theorem_}
Let $(M,\omega, \theta)$ be a compact strict LCK manifold,
$\dim_\C M=n$.
Assume that \ref{_pseudo_ef_conjecture_} holds for $M$,
that is, that the Bott--Chern class of $d^c \theta$ is
pseudo-effective. Then $M$ does not admit a balanced metric.

\hfill

\proof
Let $\omega_1$ be a balanced Hermitian form on $M$.
Then $\omega^{n-1}_1$ is closed. Therefore, \ref{_pseudo_ef_conjecture_}
implies that $\int_M \omega^{n-1}_1\wedge d^c \theta>0$.
This is impossible, because $d^c\theta = - d\theta^c$ is exact.
\endproof

\hfill

Together with \ref{_pseff_conj_true_for_known_classes_Theorem_},
this result immediately implies:

\hfill

\corollary\label{_no_balanced_on_lck_Corollary}
Let $(M,\omega, \theta)$ be a compact  non-K\"ahler LCK-manifold which 
is bimeromorphic to either  an LCK manifold with potential,  
an Oeljeklaus--Toma manifold or a Kato
manifold. Then  $(M,\omega, \theta)$ does not admit a balanced metric.

\proof \ref{_pseudo_ef_conjecture_} is true for these three classes
of LCK manifolds by \ref{_pseff_conj_true_for_known_classes_Theorem_}.
\endproof

\hfill

\remark All known examples of LCK manifolds fall in one of
the three classes described in Sections 3.2, 3.3, 3.4 (LCK
with potential, Oeljeklaus--Toma, and Kato)  or are
blow-ups of manifolds in these classes. Consequently, 
\ref{_no_balanced_on_lck_Corollary}
is at the moment the best result in this area.

\hfill

\noindent{\bf Acknowledgment:} We thank Victor Vuletescu for a careful reading of a first draft of the paper.

\hfill

{\scriptsize

	\noindent {\sc Liviu Ornea\\
		{\sc University of Bucharest, Faculty of Mathematics and Informatics, \\14
			Academiei str., 70109 Bucharest, Romania}, \\
			also:\\
		Institute of Mathematics ``Simion Stoilow" of the Romanian
		Academy\\
		21, Calea Grivitei Str.
		010702-Bucharest, Romania}\\
	{\tt lornea@fmi.unibuc.ro,   liviu.ornea@imar.ro}
	
	\hfill

	\noindent
	{\sc Misha Verbitsky\\
		{\sc Instituto Nacional de Matem\'atica Pura e
			Aplicada (IMPA) \\ Estrada Dona Castorina, 110\\
			Jardim Bot\^anico, CEP 22460-320\\
			Rio de Janeiro, RJ - Brasil }\\
{\tt verbit@impa.br}
}}
\end{document}